\documentclass[12pt,leqno,twoside]{article}

\usepackage{amssymb}
\usepackage{amsmath}
\usepackage[matrix,arrow,curve]{xy}
\renewcommand{\mod}{\operatorname{mod}}
\newcommand{\umod}{\operatorname{\underline{mod}}}
\newcommand{\Hom}{\operatorname{Hom}}
\newcommand{\Triv}{\operatorname{T}}
\newcommand{\soc}{\operatorname{soc}}
\newcommand{\charact}{\operatorname{char}}

\usepackage{amsthm}

\title{\Large\bf Derived equivalence classification of symmetric algebras of domestic type}
\author{\normalsize By Thorsten \textsc{Holm} and Andrzej
\textsc{Skowro\'nski}}
\date{}

\makeatletter

\renewcommand{\section}{\@startsection{section}{0}{\parindent}%
  {\baselineskip}{3ex}{\bf}}

\renewcommand*{\@seccntformat}[1]{%
  \csname the#1\endcsname.\quad
}

\newtheoremstyle{theorem}{}{}{\itshape}{\parindent}{\scshape}{.}{ }{}
\newtheoremstyle{example}{}{}{\upshape}{\parindent}{\scshape}{.}{ }{}

\theoremstyle{theorem}
\newtheorem{theorem}{Theorem}[section]

\newtheorem{proposition}[theorem]{Proposition}

\theoremstyle{example}

\renewenvironment{proof}{
  \trivlist
  \item[\hskip\labelsep{\scshape\hspace*{\parindent}Proof.}]
}{
  \qed
  \endtrivlist
}

\newenvironment{arablist}{
  \begin{list}{\textup{(\arabic{enumi})}}{
    \usecounter{enumi}
    \setlength{\topsep}{0pt}
    \setlength{\parsep}{0pt}
    \setlength{\itemsep}{0pt}
    \setlength{\leftmargin}{0pt}
    \settowidth{\labelwidth}{(m)}
    \setlength{\itemindent}{\parindent}
    \addtolength{\itemindent}{\labelwidth}
    \addtolength{\itemindent}{\labelsep}
  }
}{
  \end{list}
}

  \renewenvironment{abstract}{%
         \quotation
     \scriptsize
     {\bfseries \abstractname.}%
    }
        \endquotation

\makeatletter

\def\ps@paper{
  \let\@mkboth\@gobbletwo
  \def\@evenhead{%
     \parbox{\textwidth}{
       \upshape\footnotesize\thepage\hfill%
    \textsc{T. Holm} and \textsc{A. Skowro\'nski}%
       \hfill\rule[-3pt]{0pt}{15pt}}
     }
  \def\@oddhead{%
    \parbox{\textwidth}{
    \upshape\footnotesize\hfill%
    \textit{Derived equivalence classification of symmetric algebras of domestic type}%
    \hfill\thepage\rule[-3pt]{0pt}{15pt}}
  }
  \def\@oddfoot{}
  \let\@evenfoot\@oddfoot
}

\def\ps@first{
  \let\@mkboth\@gobbletwo
  \def\@evenhead{}%
  \def\@oddhead{}%
  \def\@oddfoot{\hfill\footnotesize\thepage\hfill}%
  \let\@evenfoot\@oddfoot%
}

\makeatother

\begin{document}

\maketitle

\thispagestyle{empty}

\renewcommand{\thefootnote}{}
\footnote{The second named author acknowledges supported
    from the Polish Scientific Grant KBN No. 1 P03A 018 27.}

\begin{abstract}
We give a complete derived equivalence classification of all
symmetric algebras of domestic representation type over an
algebraically closed field.
This completes previous work by R.\,Bocian and the authors,
where in this paper we solve the crucial problem of distinguishing
standard and nonstandard algebras up to derived equivalence.
Our main tool are generalized Reynolds ideals, introduced
by B.\,K\"ulshammer  for symmetric algebras in positive
characteristic, and recently shown by
A.\,Zimmermann to be invariants under derived equivalences.
\end{abstract}


\section{Introduction.}

Throughout the paper $K$ will denote a fixed algebraically closed
field. By an algebra we mean a finite dimensional $K$-algebra. For
an algebra $A$, we denote by $\mod A$ the category of finite
dimensional right $A$-modules and by $D$ the standard duality
$\Hom_K(-,K)$ on $\mod A$. An algebra $A$ is called
\textit{selfinjective} if $A \cong D(A)$ in $\mod A$, that is the
projective $A$-modules are injective. Further, an algebra $A$ is
called \textit{symmetric} if $A$ and $D(A)$ are isomorphic as
$A$-$A$-bimodules. Recall also that an algebra $A$ is symmetric if
and only if there exists an associative, symmetric, nondegenerate
$K$-bilinear form $(-,-) : A \times A \rightarrow K$. The
classical examples of selfinjective algebras (respectively,
symmetric algebras) are provided by the finite dimensional Hopf
algebras (respectively, the group algebras of finite groups).
Moreover, for any algebra $B$, the trivial extension
$\Triv(B)=B\ltimes D(B)$ of $B$ by the $B$-$B$-bimodule $D(B)$ is
a symmetric algebra, and $B$ is a factor algebra of $\Triv(B)$.
If $A$ is a selfinjective algebra, then the left and the right
socle of $A$ concide, and we denote them by $\soc(A)$. Two
selfinjective algebras $A$ and $\Lambda$ are said to be
\textit{socle equivalent} if the factor algebras $A / \soc(A)$ and
$\Lambda / \soc(\Lambda)$ are isomorphic.
For an algebra $A$, we denote by $D^b(\mod A)$ the derived
category of bounded complexes from $\mod A$. Finally, two algebras
$A$ and $\Lambda$ are said to be \textit{derived equivalent} if
the derived categories $D^b(\mod A)$ and $D^b(\mod \Lambda)$ are
equivalent as triangulated categories.

Since Happel's work \cite{Ha1} interpreting tilting theory in
terms of equivalences of derived categories, the machinery of
derived categories has been of interest to representation
theorists. In \cite{Ric1} J.\,Rickard proved his celebrated
criterion: two algebras $A$ and $\Lambda$ are derived equivalent
if and only if $\Lambda$ is the endomorphism algebra of a tilting
complex over $A$. Since a lot of interesting properties are
preserved by derived equivalences, it is for many purposes
important to classify classes of algebras  up to derived
equivalence, instead of Morita equivalence. For instance, for
selfinjective algebras the representation type is an invariant of
the derived category. Further, derived equivalent selfinjective
algebras are stably equivalent \cite{Ric2}, and hence have
isomorphic stable Auslander-Reiten quivers. It has been also
proved in \cite{Ric3} that the class of symmetric algebras is
closed under derived equivalences. Finally, we note that derived
equivalent algebras have the same number of pairwise nonisomorphic
simple modules and isomorphic centers.

One central problem of modern representation theory is the
determination of the derived equivalence classes of selfinjective
algebras of tame representation type. Recall that for a tame
algebra the indecomposable modules occur, in each dimension $d$,
in a finite number of discrete and a finite number of
one-parameter families. A distinguished class of tame algebras is
formed by the representation-finite algebras for which there are
only finitely many isomorphism classes of indecomposable modules.
In \cite{Ric2} J.\,Rickard classified the Brauer tree algebras
(for instance, representation-finite blocks of group algebras) up
to derived equivalence in connection with Brou\'e's conjecture
\cite{Br}. The derived equivalence classification of all
representation-finite selfinjective algebras has been established
by H.\,Asashiba \cite{A}. We refer also to \cite{Ho} for the
derived equivalence classification of algebras of the dihedral,
semidihedral and quaternion type (for instance,
representation-infinite tame blocks of group algebras), which are
tame and symmetric.

In this paper, we are concerned with the problem of derived
equivalence classification of all tame selfinjective algebras of
domestic representation type. Recall that for algebras of domestic
type there exists a common bound (independent of the fixed
dimension) for the numbers of one-parameter families of
indecomposable modules. The Morita equivalence classification of
these algebras splits into two cases: the standard algebras, whose
basic algebras admit simply connected Galois coverings, and the
remaining nonstandard algebras. By general theory (see \cite{BLR},
\cite{HW}, \cite{S1}, \cite{S2}), the connected standard
representation-finite (respectively, representation-infinite
domestic) selfinjective algebras are Morita equivalent to the
orbit algebras $\widehat{B}/G$ of the repetitive algebras
$\widehat{B}$ of tilted algebras $B$ of Dynkin (respectively,
Euclidean) type with respect to actions of admissible infinite
cyclic groups $G$ of automorphisms of $\widehat{B}$. The
nonstandard selfinjective algebras of domestic type are very
exceptional and are Morita equivalent to socle and geometric
deformations of the corresponding standard selfinjective algebras
of domestic type (see \cite{BS3}, \cite{Rie}, \cite{S2},
\cite{S3}, \cite{W}).

The aim of this paper is to give a complete derived equivalence
classification of all connected representation-infinite symmetric
algebras of domestic type. The Morita equivalence classification
of these algebras has been established recently in \cite{BS1},
\cite{BS2}, \cite{BS3}, \cite{S3}. In Section~\ref{sec:2} we
define (by quivers and relations) the following families of
representation-infinite domestic symmetric algebras:

\smallskip

$A(p,q)$, where $1 \leq p \leq q$,

$\Lambda(m)$, where $m \geq 2$,

$\Gamma(n)$, where $n\geq 1$,

$T(p,q)$, where $1 \leq p \leq q$,

$T(2,2,r)^{*}$, where $r \geq 2$,

$T(3,3,3)$, $T(2,4,4)$, and $T(2,3,6)$,

$\Omega(n)$, where $n \geq 1$ and $\charact K = 2$.

\smallskip

The following theorem is the main result of the paper.

\begin{theorem}
\label{th:1.1}%
The algebras $A(p,q)$, $\Lambda(m)$, $\Gamma(n)$, $T(p,q)$,
$T(2,2,r)^{*}$, $T(3,3,3)$, $T(2,4,4)$, $T(2,3,6)$, and
$\Omega(n)$ ($\charact K = 2$) form a complete set of
representatives of pairwise different derived equivalence classes
of connected representation-infinite symmetric algebras of
domestic type.
\end{theorem}

The derived equivalence classification of the standard
(respectively, nonstandard) representation-infinite symmetric
algebras of domestic type has been established in our joint papers
with R.\,Bocian \cite{BHS1} (respectively, \cite{BHS2}).
However, it remained open in these papers whether a standard
and a nonstandard algebra can be derived equivalent, or not.
(Recall that in general it is a notoriously difficult problem to
distinguish algebras up to derived equivalence. The main problem
is usually to find suitable derived invariants which are possible to
compute.)
In this paper we solve this problem, thereby completing
the derived equivalence classification of symmetric
algebras of domestic representation type.
More precisely, we prove in
Section~\ref{sec:3} that the derived equivalence classes of the
connected standard and nonstandard representation-infinite
symmetric algebras of domestic type are disjoint.
The crucial tool for proving this are the socalled
generalized Reynolds ideals defined by B.\,K\"ulshammer
in \cite{Ku1} for symmetric algebras in positive
characteristic. These sequences of ideals of the center of the
algebra have recently been
shown by A.\,Zimmermann to be invariant under derived equivalences
\cite{Z}.
This invariant is suitable for our purposes since the nonstandard symmetric
algebras of domestic type occur only in characteristic $2$. In the final
Section~\ref{sec:4} we present (for completeness) the derived
equivalence classification of all representation-finite symmetric
algebras from \cite{A}, and give an alternative proof
of the important step in Asashiba's
classification that the derived equivalence classes of the
connected standard and nonstandard representation-finite symmetric
algebras are disjoint. Our argument, using the above Reynolds ideals,
considerably simplifies the original proof in \cite{A}.

For basic background on the representation theory applied here we
refer to the books \cite{ASS}, \cite{ARS}, \cite{Ha2}, \cite{KZ}
and to the survey articles \cite{S2}, \cite{Y}.


\section{Derived normal forms of domestic symmetric algebras.}
\label{sec:2}

In order to define derived normal forms of the connected
representation-infinite domestic symmetric algebras consider the
following families of quivers
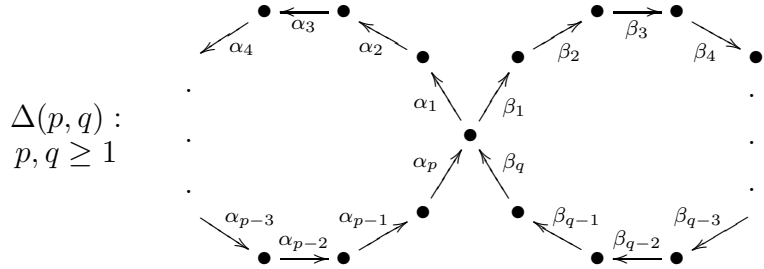
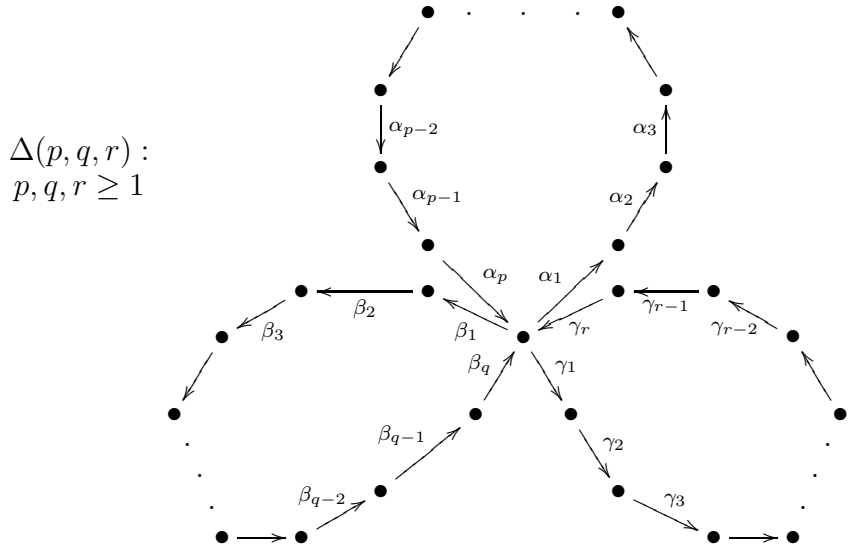
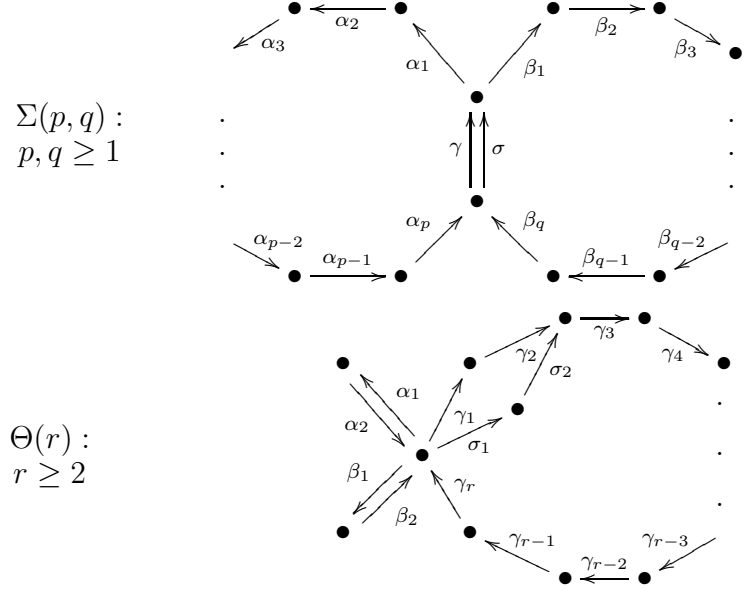
\[
\xymatrix@C=.5pc@R=.5pc{
    &&&
     \bullet \ar^{\alpha_4}[lld] && \bullet \ar^{\alpha_3}[ll] &&& &&& \bullet \ar_{\beta_3}[rr] && \bullet \ar_{\beta_4}[rrd] \\
    &
  \ar@{*{}*{\rule[2.5mm]{0mm}{0mm}.\rule[-2.5mm]{0mm}{0mm}}*{}}[dddd]
      &&&&&&
      \bullet \ar^{\alpha_2}[ull] && \bullet \ar_{\beta_2}[urr] &&&&&& \bullet
  \ar@{*{}*{\rule[2.5mm]{0mm}{0mm}.\rule[-2.5mm]{0mm}{0mm}}*{}}[dddd]
      \\
    &
    &&&&&&&&&&&&&&
    &
    \\
     \save[] +<-3pc,0pt> *\txt<4pc>{$\Delta(p,q):$\\ $p,q \geq 1$} \restore
      &&&&&&&& \bullet \ar^{\alpha_1}[luu] \ar_{\beta_1}[ruu] 
      \\
     &&&&&&&&&&&&&&&&
     \\
    &
      \ar^{\alpha_{p-3}}[rrd] &&&&&&
        \bullet \ar^{\alpha_{p}}[ruu] && \bullet \ar_{\beta_{q}}[luu]
        &&&&&& \ar_{\beta_{q-3}}[lld] \\
    &&& \bullet \ar^{\alpha_{p-2}}[rr] && \bullet \ar^{\alpha_{p-1}}[rru]
      &&& &&&
      \bullet \ar_{\beta_{q-1}}[llu] && \bullet \ar_{\beta_{q-2}}[ll] \\
}
\]
\[
\xymatrix@C=.5pc@R=.5pc{
    &&&&&& \bullet \ar[ldd] \ar@{*{ }*{\ \ .\ \ }*{ }}[rrrr] && && \bullet \\
    \\
    &&&&& \bullet \ar^{\alpha_{p-2}}[dd] &&& &&& \bullet \ar[luu] \\
    \\
\save[] +<-3pc,0pt> *\txt<4.5pc>{$\Delta(p,q,r):$\\ $p,q,r \geq
1$} \restore
    &&&&& \bullet \ar^{\alpha_{p-1}}[rdd] &&& &&& \bullet \ar^{\alpha_{3}}[uu] \\
    \\
    &&&&&& \bullet \ar^{\alpha_{p}}[rrdd] && && \bullet \ar^{\alpha_{2}}[ruu] \\
    &&& \bullet \ar^{\beta_{3}}[lld] &&& \bullet \ar^{\beta_{2}}[lll] && &&
         \bullet \ar^{\gamma_{r}}[lld] && \bullet \ar^{\gamma_{r-1}}[ll] \\
    & \bullet \ar[ldd] &&&&&&& \bullet \ar^{\alpha_{1}}[rruu] \ar^{\beta_{1}}[llu] \ar^{\gamma_{1}}[rdd] &&&&&& \bullet \ar^{\gamma_{r-2}}[llu] \\
    \\
    \bullet \ar@{*{ }*{\rule[-2mm]{0mm}{0mm}.\rule[2mm]{0mm}{0mm}}*{ }}[rddd]
        &&&&&&& \bullet \ar^{\beta_{q}}[ruu] & & \bullet \ar^{\gamma_{2}}[rdd] &&&&&&
        \bullet \ar[luu] \ar@{*{ }*{\rule[-2mm]{0mm}{0mm}.\rule[2mm]{0mm}{0mm}}*{ }}[lddd] \\
    \\
    &&&&& \bullet \ar^{\beta_{q-1}}[rruu] &&& && \bullet \ar^{\gamma_{3}}[rrd] \\
    & \bullet \ar[rr] && \bullet \ar^{\beta_{q-2}}[rru] &&&&& &&&& \bullet \ar[rr] && \bullet \\
}
\]
\[
\xymatrix@C=.45pc@R=.45pc{
    &&&
     \bullet \ar^{\alpha_3}[lld] &&& \bullet \ar^{\alpha_2}[lll] && && \bullet \ar_{\beta_2}[rrr] &&& \bullet \ar_{\beta_3}[rrd] \\
    &
      &&&&&&&
      &&&&&&& \bullet\\
     &
\ar@{*{}*{\rule[1.5mm]{0mm}{0mm}.\rule[-1.5mm]{0mm}{0mm}}*{}}[ddd]
     &&&&&&& \bullet \ar^ {\alpha_1}[uull] \ar_{\beta_1}[uurr]
&&&&&&&
\ar@{*{}*{\rule[1.5mm]{0mm}{0mm}.\rule[-1.5mm]{0mm}{0mm}}*{}}[ddd]
     \\
     \save[] +<-4pc,0pt> *\txt<4pc>{$\Sigma(p,q):$\\ $p,q \geq 1$} \restore
      \\
      \\
       &&&&&&&& \bullet \ar@<.5ex>^{\gamma}[uuu] \ar@<-.5ex>_{\sigma}[uuu]
&&&&&&&
       \\
    &
      \ar^{\alpha_{p-2}}[rrd] &&&&&&&
        &&&&&&& \ar_{\beta_{q-2}}[lld] \\
    &&& \bullet \ar^{\alpha_{p-1}}[rrr] &&& \bullet \ar^{\alpha_{p}}[rruu]
      && &&
      \bullet \ar_{\beta_{q}}[lluu] &&& \bullet \ar_{\beta_{q-1}}[lll] \\
}
\]
\[
\xymatrix@C=.5pc@R=.5pc{
    &&&&&&&& &&& \bullet \ar_{\gamma_3}[rr] && \bullet \ar_{\gamma_4}[rrd] \\
    &&&&&& \bullet \ar@<-.5ex>_{\alpha_2}[rrdd] && &
      \bullet \ar_{\gamma_2}[urr] &&&&&& \bullet
    \ar@{*{}*{\rule[3mm]{0mm}{0mm}.\rule[-3mm]{0mm}{0mm}}*{}}[dddd]
      \\
    &&&&&&&& && \bullet \ar_{\sigma_2}[ruu] &&&&&
    &
    \\
     \save[] +<-3pc,0pt> *\txt<4pc>{$\Theta(r):$\\ $r \geq 2$} \restore
      &&&&&&&& \bullet \ar@<-.5ex>_{\alpha_1}[lluu] \ar@<-.5ex>_{\beta_1}[lldd]
           \ar_{\gamma_1}[ruu] \ar_{\sigma_1}[rru] 
      \\
     &&&&&&&&&&&&&&&&
     \\
    &&&&&& \bullet \ar@<-.5ex>_{\beta_2}[rruu] && & \bullet \ar_{\gamma_{r}}[luu]
        &&&&&& \ar_{\gamma_{r-3}}[lld] \\
    &&&&&&&& &&&
      \bullet \ar_{\gamma_{r-1}}[llu] && \bullet \ar_{\gamma_{r-2}}[ll] \\
}
\]

\medskip

{\sc The algebras $A(p,q)$.}
For $1 \leq p \leq q$, denote by $A(p,q)$ the algebra given by the
quiver $\Delta(p,q)$ and the relations:
\[
 \begin{array}{l}
 \alpha_{1}\alpha_{2}\ldots\alpha_{p}\beta_{1}\beta_{2}\ldots\beta_{q}=
 \beta_{1}\beta_{2}\ldots\beta_{q}\alpha_{1}\alpha_{2}\ldots\alpha_{p},\\
 \alpha_{p}\alpha_{1}=0,\; \beta_{q}\beta_{1}=0,\\
 \alpha_{i}\alpha_{i+1}\ldots\alpha_{p}\beta_{1}\ldots\beta_{q}\alpha_{1}\ldots\alpha_{i-1}\alpha_{i}=0,\;
 2 \leq i \leq p - 1,\\
 \beta_{j}\beta_{j+1}\ldots\beta_{q}\alpha_{1}\ldots\alpha_{p}\beta_{1}\ldots\beta_{j-1}\beta_{j}=0,\;
 2\leq j \leq q - 1.
 \end{array}
\]
We note that $A(p,q)$ is a standard one-parametric symmetric
algebra of Euclidean type $\widetilde{\mathbb{A}}_{2(p+q)-3}$ (see
\cite[(5.3)(1)]{BHS1}).
\smallskip

{\sc The algebras $\Lambda(m)$.}
For $m \geq 2$, denote by $\Lambda(m)$ the algebra given by the
quiver $\Delta(1,m)$ and the relations:
\[
 \begin{array}{l}
 \alpha_{1}^2 = \left(\beta_{1}\beta_{2}\ldots\beta_{m}\right)^2,\;
     \alpha_1 \beta_{1} = 0,\;  \beta_{m} \alpha_{1} = 0,\\
 \beta_{j}\beta_{j+1}\ldots\beta_{m}\beta_{1}\ldots\beta_{m}\beta_{1}\ldots\beta_{j-1}\beta_{j}=0,\;
 2\leq j \leq m - 1.
 \end{array}
\]
We note that $\Lambda(m)$ is a standard one-parametric symmetric
algebra of Euclidean type $\widetilde{\mathbb{A}}_{2m-1}$ (see
\cite[(5.3)(2)]{BHS1}).
\smallskip

{\sc The algebras $\Gamma(n)$.}
For $n \geq 1$, denote by $\Gamma(n)$ the algebra given by the
quiver $\Delta(2,2,n)$ and the relations:
\[
 \begin{array}{l}
 \alpha_{1}\alpha_{2}=\left(\gamma_{1}\gamma_{2}\ldots\gamma_{n}\right)^2=\beta_{1}\beta_{2},\\
 \alpha_{2}\gamma_{1}=0,\; \beta_{2}\gamma_{1}=0,\; \gamma_{n}\alpha_{1}=0,\\
 \gamma_{n}\beta_{1}=0,\; \alpha_{2}\beta_{1}=0,\; \beta_{2}\alpha_{1}=0,\\
 \gamma_{j}\gamma_{j+1}\ldots\gamma_{n}\gamma_{1}\ldots\gamma_{n}\gamma_{1}\ldots\gamma_{j-1}\gamma_{j}=0,\;
 2\leq j \leq n-1.
 \end{array}
\]
We note that $\Gamma(n)$ is a standard one-parametric symmetric
algebra of Euclidean type $\widetilde{\mathbb{D}}_{2n+3}$ (see
\cite[(5.3)(3)]{BHS1}).
\smallskip

{\sc The algebras $T(p,q,r)$.}
For $2 \leq p \leq q \leq r$, denote by $T(p,q,r)$ the algebra
given by the quiver $\Delta(p,q,r)$ and the relations:
\[
 \begin{array}{l}
 \alpha_{1}\alpha_{2}\ldots\alpha_{p} = \beta_{1}\beta_{2}\ldots\beta_{q} = \gamma_{1}\gamma_{2}\ldots\gamma_{r},\\
 \beta_{q}\alpha_{1}=0,\; \gamma_{r}\alpha_{1}=0,\; \alpha_{p}\beta_{1}=0,\\
 \gamma_{r}\beta_{1}=0,\; \alpha_{p}\gamma_{1}=0,\; \beta_{q}\gamma_{1}=0,\\
 \alpha_{i}\alpha_{i+1}\ldots\alpha_{p}\alpha_{1}\ldots\alpha_{i-1}\alpha_{i}=0,\;
 2\leq i \leq p-1,\\
 \beta_{j}\beta_{j+1}\ldots\beta_{q}\beta_{1}\ldots\beta_{j-1}\beta_{j}=0,\;
 2\leq j \leq q-1,\\
 \gamma_{k}\gamma_{k+1}\ldots\gamma_{r}\gamma_{1}\ldots\gamma_{k-1}\gamma_{k}=0,\;
 2\leq k \leq r-1.
 \end{array}
\]
Observe that $T(p,q,r)$ is isomorphic to the trivial extension
algebra $\Triv(H(p,q,r))$ of the path algebra $H(p,q,r)$ of the
quiver $\Delta^{*}(p,q,r)$ obtained from $\Delta(p,q,r)$ by
deleting the arrows $\alpha_1$, $\beta_1$, $\gamma_1$. Further,
$1/p + 1/q + 1/r > 1$ if and only if $(p,q,r) = (2,2,r), (2,3,3),
(2,3,4), (2,3,5)$, or equivalently $\Delta(p,q,r)$ is a Dynkin
quiver of type $\mathbb{D}_{r+2}$, $\mathbb{E}_{6}$,
$\mathbb{E}_{7}$, or $\mathbb{E}_{8}$, respectively. In this case,
$T(p,q,r)$ is a standard representation-finite symmetric algebra
(see \cite{HW}). Similarly, $1/p + 1/q + 1/r = 1$ if and only if
$(p,q,r) = (3,3,3), (2,4,4), (2,3,6)$, or equivalently
$\Delta(p,q,r)$ is an Euclidean quiver of type
$\widetilde{\mathbb{E}}_{6}$, $\widetilde{\mathbb{E}}_{7}$, or
$\widetilde{\mathbb{E}}_{8}$, respectively. In this case,
$T(p,q,r)$ is a standard $2$-parametric symmetric algebra (see
\cite{ANS}).
\smallskip

{\sc The algebras $T(p,q)$.}
For $1 \leq p \leq q$, denote by $T(p,q)$ the algebra given by the
quiver $\Sigma(p,q)$ and the relations:
\[
 \begin{array}{l}
 \alpha_{1}\alpha_{2}\ldots\alpha_{p}\gamma =
    \beta_{1}\beta_{2}\ldots\beta_{q}\sigma,\\
 \gamma\alpha_{1}\alpha_{2}\ldots\alpha_{p} =
    \sigma\beta_{1}\beta_{2}\ldots\beta_{q},\\
 \alpha_{p}\sigma=0,\; \sigma\alpha_{1}=0,\; \beta_{q}\gamma=0,\;
\gamma\beta_{1}=0,\\
 \alpha_{i}\alpha_{i+1}\ldots\alpha_{p}\gamma\alpha_{1}\ldots
\alpha_{i-1}\alpha_{i}=0,\;
 2\leq i \leq p-1,\\
 \beta_{j}\beta_{j+1}\ldots\beta_{q}\sigma\beta_{1}\ldots\beta_{j-1}
\beta_{j}=0,\;
 2\leq j \leq q-1.\\
 \end{array}
\]
Then $T(p,q)$ is isomorphic to the trivial extension algebra
$\Triv(H(p,q))$ of the path algebra $H(p,q)$ of the quiver
$\Sigma^{*}(p,q)$ of Euclidean type
$\widetilde{\mathbb{A}}_{p+q-1}$ obtained from $\Sigma(p,q)$ by
deleting the arrows $\gamma$ and $\sigma$. In particular, $T(p,q)$
is a standard $2$-parametric symmetric algebra (see \cite{AS}).
\smallskip

{\sc The algebras $T(2,2,r)^{*}$.}
For $r \geq 2$, denote by $T(2,2,r)^{*}$ the algebra given by the
quiver $\Theta(r)$ and the relations:
\[
 \begin{array}{l}
 \alpha_{1}\alpha_{2} = \beta_{1}\beta_{2} = \gamma_{1}\gamma_{2}\ldots
\gamma_{r}
    ,\; \gamma_1\gamma_2 = \sigma_1\sigma_2,\\
 \gamma_r\alpha_{1}=0,\; \beta_2\alpha_1=0,\; \gamma_r\beta_1=0,\;
\alpha_2\beta_1=0,\\
 \alpha_2\gamma_1=0,\; \alpha_2\sigma_1=0,\; \beta_2\gamma_1=0,\;
\beta_2\sigma_1=0,\\
 \alpha_2\alpha_1\alpha_2=0,\; \beta_2\beta_1\beta_2 = 0,\\
 \gamma_{2}\gamma_{3}\ldots\gamma_{r}\sigma_1=0,\; \sigma_2\gamma_{3}
\ldots\gamma_{r}\gamma_{1}=0,\\
 \gamma_{k}\gamma_{k+1}\ldots\gamma_{r}\gamma_{1}\gamma_{2}\ldots
\gamma_{k-1}\gamma_{k}=0,\;
 3\leq k \leq r-1.
 \end{array}
\]
Then $T(2,2,r)^{*}$ is isomorphic to the trivial extension algebra
$\Triv(H(2,2,r)^{*})$ of the path algebra $H(2,2,r)^{*}$ of the
quiver $\Theta^{*}(r)$ of Euclidean type
$\widetilde{\mathbb{D}}_{r+1}$ obtained from $\Theta(r)$ by
deleting the arrows $\alpha_1$, $\beta_1$, $\gamma_1$ and
$\sigma_1$. In particular, $T(2,2,r)^{*}$ is a standard
$2$-parametric symmetric algebra (see \cite{ANS}).
\smallskip

{\sc The algebras $\Omega(n)$.}
For $n \geq 1$, denote by $\Omega(n)$ the algebra given by the
quiver $\Delta(1,n)$ and the relations:
\[
 \begin{array}{l}
 \alpha_{1}\beta_{1}\beta_{2}\ldots\beta_{n} + \beta_{1}\beta_{2}\ldots
\beta_{n}\alpha_{1} = 0,\\
 \alpha_{1}^2 = \alpha_{1}\beta_{1}\beta_{2}\ldots\beta_{n},\;
\beta_{n}\beta_{1} = 0,\\
 \beta_{j}\beta_{j+1}\ldots\beta_{n}\alpha_{1}\beta_{1}\beta_{2}\ldots
\beta_{j-1}\beta_{j}=0,\;
 2\leq j \leq n-1.\\
 \end{array}
\]
Then $\Omega(n)$ is a nonstandard one-parametric selfinjective
algebra. Furthermore, $\Omega(n)$ is a symmetric algebra if and
only if $\charact K = 2$. Moreover, in this symmetric case,
$\Omega(n)$ is socle equivalent to the algebra $\Omega(n)^{\prime}
= A(1,n)$, called the standard form of $\Omega(n)$ (see
\cite{BS3}).

The following derived equivalence classifications of standard
symmetric algebras of domestic type has been established in
\cite[Theorems 1 and 2]{BHS1}.

\begin{theorem}
\label{th:2.1}%
The algebras $A(p,q)$, $\Lambda(m)$, $\Gamma(n)$, $T(p,q)$,
$T(2,2,r)^{*}$, $T(3,3,3)$, $T(2,4,4)$, $T(2,3,6)$ form a complete
set of representatives of pairwise different derived equivalence
classes of connected, standard, representation-infinite symmetric
algebras of domestic type.
\end{theorem}

It has been proved in \cite{S3} (see also \cite{BS3}) that the
nonstandard representation-infinite symmetric algebras of domestic
type occur only in characteristic $2$. Furthermore, the following
theorem, proved in \cite[Theorem~1]{BHS2}, gives the derived
equivalence classification of these algebras.

\begin{theorem}
\label{th:2.2}%
The algebras $\Omega(n)$, for $n \geq 1$ and $\charact K = 2$,
form a complete set of representatives of pairwise different
derived equivalence classes of connected, nonstandard,
representation-infinite symmetric algebras of domestic type.
\end{theorem}


\section{Generalized Reynolds ideals.}
\label{sec:Reynolds}

In this section we briefly recall the definition of the sequence
of generalized Reynolds ideals. For more details on this invariant
we refer to \cite{Ku1}, \cite{Ku2}, \cite{HHKM}, \cite{Z}.

Let $K$ an algebraically closed field of characteristic $p > 0$.
Let $A$ be a finite dimensional symmetric $K$-algebra with
associative, symmetric, nondegenerate $K$-bilinear form $(-,-) : A
\times A \rightarrow K$. For a $K$-subspace $M$ of $A$, denote by
$M^{\bot}$ the orthogonal complement of $M$ inside $A$ with
respect to the form $(-,-)$. Moreover, let $K(A)$ be the $K$-subspace
of $A$ generated by all commutators $[a,b]:=a b - b a$, for any
$a, b \in
A$. Then for any $n \geq 0$ set
$$T_n (A) = \left\{ x \in A \mid x^{p^n} \in K(A)\right\}.$$
Then, by \cite{Ku1}, the orthogonal
complements $T_n (A)^{\bot}, n \geq 0$, are ideals of the center
$Z(A)$ of $A$, called generalized Reynolds ideals.
They form a descending sequence
$$Z(A) = T_0(A)^{\perp} \supseteq T_1(A)^{\perp} \supseteq T_2(A)^{\perp}
\supseteq T_3(A)^{\perp} \supseteq \ldots$$
In fact, B.\,K\"ulshammer proved in \cite{Ku2} that the
equation $(\xi_n(z),x)^{p^n} = (z,x^{p^n})$ for any $x, z \in
Z(A)$ defines a mapping $\xi_n : Z(A) \rightarrow Z(A)$ such that
$\xi_n(A) = T_n (A)^{\bot}$.

Then we have the following theorem proved recently by
A.\,Zimmermann \cite{Z}.

\begin{proposition}
\label{prop:zimmermann} %
Let $A$ and $B$ be derived equivalent symmetric algebras over an
algebraically closed field of positive characteristic $p$. Then
there is an isomorphism $\varphi : Z(A) \rightarrow Z(B)$ of the
centers of $A$ and $B$ such that $\varphi(T_n (A)^{\bot}) = T_n
(B)^{\bot}$ for all positive integers $n$.
\end{proposition}

\begin{proof}
See \cite[Theorem~1]{Z}.
\end{proof}

Hence the sequence of generalized Reynolds ideals gives a new
derived invariant, for symmetric algebras over algebraically
closed fields of positive characteristic.
\smallskip

In the next section we shall use this invariant for proving our
main result. The algebras occurring in our context are all given
by a quiver with relations (bound quiver) and a basis of the
algebra is provided by the set of all pairwise distinct (modulo
the ideal generated by the imposed relations) nonzero paths of the
quiver.

The following simple observation will turn out to be useful.

\begin{proposition}
\label{prop:form} %
Let $A=KQ/I$ be a symmetric bound quiver algebra, and assume that
a $K$-basis ${\mathcal B}$ of $A$ is given by the pairwise
distinct nonzero paths of the quiver $Q$ (modulo the ideal $I$).
Then the following statements hold:
\begin{enumerate}
\item[{(1)}] An associative nondegenerate symmetric $K$-bilinear
form $(-,-)$ for $A$ is given as follows
$$(x,y) = \left\{
\begin{array}{ll} 1 & \mbox{if $xy\in\soc(A)$} \\
                  0 & \mbox{otherwise}
\end{array} \right.
$$
for $x,y\in {\mathcal B}$.
\item[{(2)}] For any $n\ge 0$, the socle $\soc(A)$ is contained in
the generalized Reynolds ideal $T_n(A)^{\perp}$.
\end{enumerate}
\end{proposition}

\begin{proof} (1)
It is well-known (see \cite[Section~2]{Y}) that an algebra $A$ is
symmetric if and only if there is a $K$-linear form $\psi : A
\rightarrow K$ such that $\psi (a b) = \psi(b a)$ for all elements
$a, b \in A$ and the kernel of $\psi$ contains no nonzero left or
right ideal of $A$. Moreover,  for such a (symmetrizing) form
$\psi : A \rightarrow K$, the $K$-bilinear form $(-,-) : A \times
A \rightarrow K$ given by $(a,b) = \psi (a b)$ for all $a,b \in A$
is an associative, symmetric, nondegenerate form. For the
symmetric algebra $A = K Q / I$ considered in the proposition, we
may take the symmetrizing form $\varphi : A \rightarrow K$ which
assigns $1$ to any nonzero residue class of a path in $Q$ in $A =
K Q / I$ from the socle $\soc (A)$, and $0$ to the residue classes
of the remaining paths of $Q$. Then the bilinear form $(-,-)$
associated to $\varphi$ satisfies the required statement (1).

\smallskip

(2) By \cite{HHKM} we have for any symmetric algebra $A$ that
$$\bigcap_{n=0}^{\infty} T_n(A)^{\perp} = \soc(A) \cap Z(A).$$
But for the algebras described in the proposition we always have
$\soc(A)\subseteq Z(A)$.
\end{proof}


\section{Proof of the main result.}
\label{sec:3}

The aim of this section is to give the proof of
Theorem~\ref{th:1.1}, applying Theorems \ref{th:2.1} and
\ref{th:2.2}. We need some general facts.

For a selfinjective algebra $A$, we denote by $\Gamma^s_A$ the
\textit{stable Auslander-Reiten quiver} of $A$, obtained from the
Auslander-Reiten quiver $\Gamma_A$ of $A$ by removing the
projective-injective vertices and the arrows attached to them.
Recall that two selfinjective algebras $A$ and $B$ are called
\textit{stably equivalent} if their stable module categories
$\umod A$ and $\umod B$ are equivalent.

\begin{proposition}
\label{prop:3.1}%
Let $A$ and $B$ be stably equivalent connected selfinjective
algebras of Loewy length at least $3$. Then the stable
Auslander-Reiten quivers $\Gamma^s_A$ and $\Gamma^s_B$ are
isomorphic.
\end{proposition}

\begin{proof}
See \cite[Corollary~X.1.9]{ARS}.
\end{proof}

\begin{proposition}
\label{prop:3.2}%
Let $A$ and $\Lambda$ be derived equivalent selfinjective
algebras. Then $A$ and $\Lambda$ are stably equivalent.
\end{proposition}

\begin{proof}
See \cite[Corollary~2.2]{Ric2}.
\end{proof}

We know from Theorem~\ref{th:2.1} (respectively,
Theorem~\ref{th:2.2}) that the  algebras $A(p,q)$, $\Lambda(m)$,
$\Gamma(n)$, $T(p,q)$, $T(2,2,r)^{*}$, $T(3,3,3)$, $T(2,4,4)$,
$T(2,3,6)$ (respectively, $\Omega(n)$, for $\charact K = 2$) form
a complete set of representatives of pairwise different derived
equivalence classes of connected standard (respectively,
nonstandard) representation-infinite symmetric algebras of
domestic type. Moreover, these algebras are basic, connected and
of Loewy length at least $3$. It also follows from Propositions
\ref{prop:3.1} and \ref{prop:3.2} that the stable Auslander-Reiten
quivers of two derived equivalent connected selfinjective algebras
of Loewy length at least $3$ are isomorphic. In order to
distinguish the derived equivalence classes of the standard and
nonstandard representation-infinite symmetric algebras of domestic
type, we need the shapes of the stable Auslander-Reiten quivers of
algebras occuring in Theorems \ref{th:2.1} and \ref{th:2.2}.

\begin{proposition}
The following statements hold:
\begin{arablist}
  \item
    $\Gamma^s_{A(p,q)}$ consists of an Euclidean component of type
    $\mathbb{Z}\widetilde{\mathbb{A}}_{2(p+q)-3}$ and a
    $\mathbb{P}_1(K)$-family of stable tubes of tubular type
    $(2p-1,2q-1)$.
  \item
    $\Gamma^s_{\Lambda(n)}$ consists of an Euclidean component of type
    $\mathbb{Z}\widetilde{\mathbb{A}}_{2n-1}$ and a
    $\mathbb{P}_1(K)$-family of stable tubes of tubular type
    $(n,n)$.
  \item
    $\Gamma^s_{\Gamma(n)}$ consists of an Euclidean component of type
    $\mathbb{Z}\widetilde{\mathbb{D}}_{2n+3}$ and a
    $\mathbb{P}_1(K)$-family of stable tubes of tubular type
    $(2,2,2n+1)$.
\end{arablist}
\end{proposition}

\begin{proof}
See \cite[Proposition~5.3]{BHS1}.
\end{proof}

\begin{proposition}
The following statements hold:
\begin{arablist}
  \item
    $\Gamma^s_{T(p,q)}$ consists of two Euclidean components of type
    $\mathbb{Z}\widetilde{\mathbb{A}}_{p+q-1}$ and two
    $\mathbb{P}_1(K)$-families of stable tubes of tubular type
    $(p,q)$.
  \item
    $\Gamma^s_{T(2,2,r)^{*}}$ consists of two Euclidean components of type
    $\mathbb{Z}\widetilde{\mathbb{D}}_{r+1}$ and two
    $\mathbb{P}_1(K)$-families of stable tubes of tubular type
    $(2,2,r+2)$.
  \item
    $\Gamma^s_{T(3,3,3)}$ consists of two Euclidean components of type
    $\mathbb{Z}\widetilde{\mathbb{E}}_{6}$ and two
    $\mathbb{P}_1(K)$-families of stable tubes of tubular type
    $(2,3,3)$.
  \item
    $\Gamma^s_{T(2,4,4)}$ consists of two Euclidean components of type
    $\mathbb{Z}\widetilde{\mathbb{E}}_{7}$ and two
    $\mathbb{P}_1(K)$-families of stable tubes of tubular type
    $(2,3,4)$.
  \item
    $\Gamma^s_{T(2,3,6)}$ consists of two Euclidean components of type
    $\mathbb{Z}\widetilde{\mathbb{E}}_{8}$ and two
    $\mathbb{P}_1(K)$-families of stable tubes of tubular type
    $(2,3,5)$.
\end{arablist}
\end{proposition}

\begin{proof}
We know that $T(p,q)$, $T(2,2,r)^{*}$, $T(3,3,3)$, $T(2,4,4)$,
$T(2,3,6)$ are the trivial extension algebras of the hereditary
algebras $H(p,q)$, $H(2,2,r)^{*}$, $H(3,3,3)$, $H(2,4,4)$,
$H(2,3,6)$ of Euclidean types $\widetilde{\mathbb{A}}_{p+q-1}$,
$\widetilde{\mathbb{D}}_{r+1}$, $\widetilde{\mathbb{E}}_{6}$,
$\widetilde{\mathbb{E}}_{7}$, $\widetilde{\mathbb{E}}_{8}$,
respectively. Then the required statements follow from the
structure of the Auslander-Reiten quivers of the hereditary
algebras of Euclidean types (see \cite[(3.6)]{Rin}) and the
description of the Auslander-Reiten quivers of the trivial
extensions of the hereditary algebras given in
\cite[Theorem~3.4]{T}, \cite[Theorem~2.5.2]{Y} (see also
\cite{AS}, \cite{ANS}).
\end{proof}

\begin{proposition}
The stable Auslander-Reiten quiver $\Gamma^s_{\Omega(n)}$ of
$\Omega(n)$ consists of an Euclidean component of type
$\mathbb{Z}\widetilde{\mathbb{A}}_{2n-1}$ and a
$\mathbb{P}_1(K)$-family of stable tubes of tubular type $(2n-1)$.
\end{proposition}

\begin{proof}
See \cite[Proposition~2.2]{BHS2}.
\end{proof}

As a direct consequence of the above three propositions we obtain
the following fact.

\begin{proposition}
Let $A$ be an algebra of one of the forms  $A(p,q)$, $\Lambda(n)$,
$\Gamma(n)$, $T(p,q)$, $T(2,2,r)^{*}$, $T(3,3,3)$, $T(2,4,4)$,
$T(2,3,6)$, and let $B$ be an algebra of the form $\Omega(n)$. Assume
that the stable Auslander-Reiten quivers of $A$ and $B$ are
isomorphic. Then $A = A(1,n)$ and $B = \Omega(n)$ for some $n \geq
1$.
\end{proposition}

Therefore the following proposition completes the proof of
Theorem~\ref{th:1.1}.

\begin{proposition}
Let $K$ be an algebraically closed field of characteristic $2$.
Then, for any $n \geq 1$, the symmetric algebras $\Omega(n)$ and
$\Omega(n)^{\prime} = A(1,n)$ are not derived equivalent.
\end{proposition}

\begin{proof}
Let us denote by $\Omega$ either of the algebras
$\Omega(n)$ or $\Omega(n)'$.

We shall compute the series of generalized Reynolds ideals for the
symmetric algebras $\Omega$,
$$Z(\Omega)\supseteq T_1(\Omega)^{\perp} \supseteq T_2(\Omega)^{\perp}
\supseteq
\ldots$$
as described in Section \ref{sec:Reynolds}.

We shall show that the ideals in these series have different dimensions
for the algebras $\Omega(n)$ and $\Omega(n)^{\prime} = A(1,n)$.
Since the series of Reynolds ideals is a derived invariant
(see Proposition \ref{prop:zimmermann}), we can then distinguish these
algebras up to derived equivalence.
\smallskip

The centers of these algebras have dimension $n+2$ as vector space
over $K$. More precisely, it is straightforward to check that a
$K$-basis is given as follows
$$Z(\Omega)=\langle 1, \beta_1\ldots\beta_n, s_1:=\alpha\beta_1\ldots\beta_n,
s_j:=\beta_j\ldots\beta_n\alpha\beta_1\ldots\beta_{j-1}~(2\le j\le n)\rangle_K
$$
where we abbreviate $\alpha = \alpha_1$.
Note that the latter $n$ elements $s_1,\ldots,s_n$
form a basis of the socle of $\Omega$.

Since we are dealing with characteristic 2, we have
$$T_1(\Omega):=\{x\in \Omega\,\mid\,x^2\in K(\Omega)\}.$$
Recall that $K(\Omega)$ is the subspace of the algebra $\Omega$
generated by all commutators $[x,y]=xy-yx$, where $x,y\in \Omega$.
Now consider the first generalized Reynolds ideal
$$T_1(\Omega)^{\perp}:=\{y\in Z(\Omega)\,\mid\,(x,y)=0~\mbox{for all
$x\in T_1(\Omega)$}\},$$
where $(-,-)$ is the nondegenerate symmetric $K$-bilinear form for the
symmetric algebra $\Omega$, as defined in Proposition \ref{prop:form}.
Note that for such a basic symmetric algebra $A$, the
socle is contained in any Reynolds ideal $T_m(A)^{\perp}$
(see Proposition~\ref{prop:form}(2)).

We consider the following sequence of ideals
$$\soc(\Omega) \subseteq T_1(\Omega)^{\perp} \subset Z(\Omega).$$
Here, the second inclusion is strict, since $1$ is not contained
in any Reynolds ideal of $\Omega$. In fact, $\soc(\Omega)\subseteq
T_1(\Omega)$, and $(1,s)=1$ for every $s\in\soc(\Omega)$.

On the other hand, the socle of $\Omega$ has only codimension 2 in the
center $Z(\Omega)$, leaving us with $\beta_1\ldots\beta_n$ as the crucial
basis element to check.

But the element $\beta_1\ldots\beta_n$ is easily seen to be
orthogonal to all basis elements in the ideal generated
by the arrows of the quiver, except to $\alpha$.
In fact, $(\alpha,\beta_1\ldots\beta_n)=1$ since
$\alpha\beta_1\ldots\beta_n$ belongs to the socle of $\Omega$.
\smallskip

Now we have to consider the algebras $\Omega(n)$ and $\Omega(n)'$
separately. Note that the distinction in the relations is that
$\alpha^2=\alpha\beta_1\ldots\beta_n$ is nonzero in $\Omega(n)$,
whereas $\alpha^2=0$ in $\Omega(n)'$.

For $\Omega(n)$,
the crucial fact to observe is that $\alpha\not\in T_1(\Omega(n))$.
In fact, $\alpha^2=\alpha\beta_1\ldots\beta_n$ is nonzero, and it
can be checked that it can not be written as a linear combination
of commutators.
But this implies that $\beta_1\ldots\beta_n\in T_1(\Omega(n))^{\perp}$.
So we get the following series of ideals and their codimensions:
$$\soc(\Omega(n)) \underbrace{\mbox{~~$\subset$~~}}_{\mbox{{\small 1}}}
T_1(\Omega(n))^{\perp} = \langle \beta_1\ldots
\beta_n, \soc(\Omega(n))\rangle_K
\underbrace{\mbox{~~$\subset$~~}}_{\mbox{{\small 1}}} Z(\Omega(n)).$$

On the other hand, for $\Omega(n)'$, we have $\alpha\in T_1(\Omega(n)')$,
since $\alpha^2=0$. Since $\beta_1\ldots\beta_n$ is not orthogonal
to $\alpha$, we conclude that $\beta_1\ldots\beta_n\not\in
T_1(\Omega(n)')^{\perp}$. Hence the corresponding series of ideals
for $\Omega(n)'$ takes the form
$$\soc(\Omega(n)') \underbrace{\mbox{~~$=$~~}}_{\mbox{{\small 0}}}
T_1(\Omega(n)')^{\perp}
\underbrace{\mbox{~~$\subset$~~}}_{\mbox{{\small 2}}} Z(\Omega(n)').$$
Since the series of generalized Reynolds ideals, and in particular
the codimensions occurring, is invariant under derived equivalences,
we can now conclude that the nonstandard algebra $\Omega(n)$ and
the standard algebra $\Omega(n)'$ are not derived equivalent.
\end{proof}


\section{Derived normal forms of representation-finite symmetric algebras.}
\label{sec:4}

In \cite[Theorem~2.2]{A} H.\,Asashiba proved that the derived
equivalence classes of connected representation-finite standard
(respectively, nonstandard) selfinjective algebras are determined
by the combinatorial data called the types, and the derived
equivalence classes of the standard and nonstandard
representation-finite selfinjctive algebras are disjoint.
Furthermore, by the results of C.\,Riedtmann
\cite[Proposition~5.7(b)]{Rie} and J.\,Waschb\"usch \cite{W}, the
nonstandard representation-finite selfinjective algebras are
symmetric, are given by some Brauer quivers, and occur only in
characteristic $2$ (see also \cite[(3.6)--(3.8)]{S3}).

The aim of this section is to present the derived equivalence
classification of all connected representation-finite symmetric
algebras by quivers and relations.

For $m,n \geq 1$, denote by $N_n^{m n}$ the algebra given by the
quiver
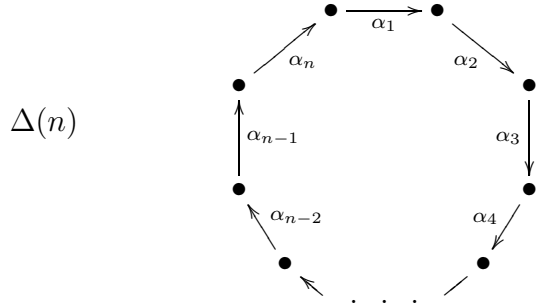
\[
\xymatrix@C=.45pc@R=.45pc{
    && && \bullet \ar_{\alpha_1}[rrr] &&& \bullet \ar_{\alpha_2}[rrdd] \\
   \\
    && \bullet \ar_{\alpha_n}[uurr]
    &&&&&&&
    \bullet \ar_{\alpha_3}[ddd]
     \\
     \save[] +<-4pc,0pt> *\txt<4pc>{$\Delta(n)$} \restore
      \\
      \\
       && \bullet \ar_{\alpha_{n-1}}[uuu] &&&&&&& \bullet \ar_{\alpha_4}[ldd]
       \\
    \\
      && &
      \bullet \ar_{\alpha_{n-2}}[luu] &&&&& \bullet \ar@{-}[ld] \\
      && && \ar[lu] &&& \ar@{{}*{\ .\ }{}}[lll] \\
}
\]
and the relations:
\[
    \left( \alpha_i \alpha_{i+1} \dots \alpha_n
    \alpha_1 \dots \alpha_{i-1}\right)^m \alpha_i = 0,\;
    1 \leq i \leq n.
\]
It is well-known that these algebras form a complete set of
representatives of the Morita equivalence classes of the symmetric
Nakayama algebras. In \cite[Theorem~4.2]{Ric2} J.\,Rickard proved
that they form a complete set of representatives of the derived
equivalence classes of the Brauer tree algebras (which occur in
the representation theory of representation-finite blocks of group
algebras). Finally, we note also that the algebra $N_n^n$ is the
trivial extension algebra $\Triv(H(n))$ of the path algebra $H(n)$
of the quiver $\Delta(n)^{*}$ obtained from $\Delta(n)$ by
deleting the arrow $\alpha_n$.

For $m \geq 2$, denote by $D(m)^{\prime}$ the algebra given by the
quiver $\Delta(1,m)$ and the relations:
\[
 \begin{array}{l}
 \alpha_{1}^2 = \beta_{1}\beta_{2}\ldots\beta_{m},\;
     \beta_{m} \beta_{1} = 0,\\
 \beta_{i}\beta_{i+1}\ldots\beta_{m}\alpha_1\beta_{1}\ldots\beta_{i-1}\beta_{i}=0,\;
 2\leq i \leq m-1.
 \end{array}
\]
Then $D(m)^{\prime}$ is a standard representation-finite symmetric
algebra of Dynkin type $\mathbb{D}_{3m}$ (see \cite{Rie}).

In Section~\ref{sec:2} we defined also the trivial extension
algebras $T(2,2,r), r \geq 2$, $T(2,3,3)$, $T(2,3,4)$, $T(2,3,5)$
of the hereditary algebras $H(2,2,r), r \geq 2$, $H(2,3,3)$,
$H(2,3,4)$, $H(2,3,5)$ of Dynkin types $\mathbb{D}_{r+2}$,
$\mathbb{E}_{6}$, $\mathbb{E}_{7}$, $\mathbb{E}_{8}$,
respectively. In \cite[Theorems 3.1 and 3.7]{HW} D.\,Hughes and
J.\,Waschb\"usch proved that the trivial extension $\Triv(B)$ of a
connected algebra $B$ is representation-finite if and only if $B$
is a tilted algebra of Dynkin type. Moreover, in
\cite[Theorem~3.1]{Ric2} J.\,Rickard proved that if $A$ and $B$
are derived equivalent algebras then their trivial extensions
$\Triv(A)$ and $\Triv(B)$ are also derived equivalent. Therefore,
the trivial extensions algebras $N_n^n = T(H(n)), n \geq 1$,
$T(2,2,r), r \geq 2$, $T(2,3,3)$, $T(2,3,4)$, $T(2,3,5)$ form a
complete set of representatives of pairwise different derived
equivalence classes of the connected representation-finite trivial
extension algebras $\Triv(B)$. We also note that by
\cite[Section~1]{BLR} all connected representation-finite
symmetric algebras of Dynkin types $\mathbb{E}_{6}$,
$\mathbb{E}_{7}$, $\mathbb{E}_{8}$ are actually the trivial
extension algebras of tilted algebras of Dynkin types
$\mathbb{E}_{6}$, $\mathbb{E}_{7}$, $\mathbb{E}_{8}$.

In \cite{Rie} C.\,Riedtmann proved that the Morita equivalence
classes of the connected standard representation-finite symmetric
algebras of Dynkin type $\mathbb{D}_n$, which are not trivial
extension algebras, are given by some (looped) Brauer trees (see
also \cite[Theorem~3.11]{S3}). Then applying Rickard's
constructions from \cite[Section~4]{Ric2} one easily proves that
the algebras $D(m)^{\prime}, m \geq 2$, give a complete set of
representatives of the derived equivalence classes of these
symmetric algebras.

Applying the combinatorial descriptions of the stable
Auslander-Reiten quivers of representation-finite selfinjective
algebras (see \cite[Section~2]{A}, \cite[Section~1]{BLR}) one
easily shows that the stable Auslander-Reiten quivers of the
standard symmetric algebras $N_n^{m n}$, $m,n\geq 1$,
$D(m)^{\prime}, m \geq 2$, $T(2,2,r), r \geq 2$, $T(2,3,3)$,
$T(2,3,4)$, $T(2,3,5)$ are pairwise nonisomorphic. Therefore,
summing up the above discussion and invoking Propositions
\ref{prop:3.1} and \ref{prop:3.2}, we obtain the following result.

\begin{proposition}
The algebras $K$, $N_n^{m n}$, $m,n\geq 1$, $D(m)^{\prime}, m \geq
2$, $T(2,2,r), r \geq 2$, $T(2,3,3)$, $T(2,3,4)$, $T(2,3,5)$ form
a complete set of representatives of pairwise different derived
equivalence classes of connected, standard, representation-finite
symmetric algebras.
\end{proposition}

For $m \geq 2$, denote by $D(m)$ the algebra given by the quiver
$\Delta(1,m)$ and the relations:
\[
 \begin{array}{l}
 \alpha_{1}^2 = \beta_{1}\beta_{2}\ldots\beta_{m},\;
     \beta_{m} \beta_{1} = \beta_{m} \alpha_{1} \beta_{1},\\
 \beta_{i}\beta_{i+1}\ldots\beta_{m}\alpha_1\beta_{1}\ldots\beta_{i-1}\beta_{i}=0,\;
 1 \leq i \leq m.
 \end{array}
\]
If $\charact K \neq 2$,
then $D(m)$ is a standard representation-finite symmetric algebra
of Dynkin type $\mathbb{D}_{3m}$ and $D(m) \cong D(m)^{\prime}$.
On the other hand, if $\charact K = 2$, then
$D(m)$ is nonstandard (see \cite[(5.7)]{Rie} or \cite{W}).
Further, $D(m)$ and $D(m)^{\prime}$ are socle equivalent, and
$D(m)^{\prime}$ is called the standard form of $D(m)$.

Applying Rickard's method from \cite[Section~4]{Ric2},
H.\,Asashiba proved in \cite[Section~7]{A} the following
proposition.

\begin{proposition}
The algebras $D(m)^{\prime}, m \geq 2$, for $\charact K = 2$, form
a complete set of representatives of pairwise different derived
equivalence classes of connected nonstandard representation-finite
symmetric algebras.
\end{proposition}

Applying again the combinatorial descriptions of the stable
Auslander-Reiten quivers of representation-finite selfinjective
algebras one easily deduces the following fact.

\begin{proposition}
Let $A$ be an algebra of one of the forms $K$, $N_n^{m n}$,
$D(m)^{\prime}$, $T(2,2,r)$, $T(2,3,3)$, $T(2,3,4)$, $T(2,3,5)$,
and $B$ an algebra of the form $D(m)$. Assume that the stable
Auslander-Reiten quivers of $A$ and $B$ are isomorphic. Then $A =
D(m)^{\prime}$ and $B = D(m)$ for some $m \geq 2$.
\end{proposition}

In \cite[Section~3]{A} H.\,Asashiba proved that, if $\charact K =
2$, then the algebras $D(m)$ and $D(m)^{\prime}$ are not stably
(hence not derived) equivalent. We will give an alternative,
simplified, proof of the latter fact, applying
Proposition~\ref{prop:zimmermann}.

\begin{proposition}
Let $K$ be an algebraically closed field of characteristic $2$.
Then, for any $m \geq 2$, the symmetric algebras $D(m)$ and
$D(m)^{\prime}$ are not derived equivalent.
\end{proposition}

\begin{proof}
For a fixed $m \geq 2$, we denote by $D$ either of the algebras
$D(m)$ or $D(m)'$.
Note that the socle $\soc (D)$ of the algebra $D$ has as $K$-basis
the following $m$ elements
$$s_1:=\alpha\beta_1\ldots\beta_n=\alpha^3=\beta_1\ldots\beta_n\alpha,$$
$$s_j:=\beta_j\ldots\beta_m\alpha\beta_1\ldots\beta_{j-1},~~
\mbox{for $j=2,\ldots,m$},$$
where we abbreviate $\alpha = \alpha_1$.

Then it is straightfoward to verify that
the center of the algebra $D$ is as $K$-vector space generated
by the following basis
$$Z(D)=\langle 1,\beta_1\ldots\beta_n,s_1,s_2,\ldots,s_m\rangle_K.$$
We shall study the series of ideals of the center
$$\soc(D)\subseteq T_1(D)^{\perp} \subseteq Z(D).$$
Since the socle only has codimension 2 in the center, and since
the unit 1 can not be contained in the Reynolds ideal
$T_1(D)^{\perp}$, the crucial question is whether
$\beta_1\ldots\beta_m$ is contained in $T_1(D)^{\perp}$, or not.

Since we are dealing with characteristic 2,
recall that
$$T_1(D):=\{x\in D\,\mid\,x^2\in K(D)\},$$
where $K(D)$ is the subspace of $D$ generated by all commutators.
\smallskip

First we consider the standard algebra $D(m)'$. We have the
relation $\beta_m\beta_1=0$, so
$$\alpha^2 = \beta_1\ldots\beta_m = [\beta_1,\beta_2\ldots\beta_m]
\in K(D(m)').$$

Thus, $\alpha\in T_1(D(m)')$. But by the definition of the
nondegenerate symmetric bilinear form on $D(m)'$, given in
Proposition~\ref{prop:form}, we have
$(\alpha,\beta_1\ldots\beta_m) =1$, because the product
$\alpha\beta_1\ldots\beta_m=s_1$ is a nonzero socle element.
Therefore, $\beta_1\ldots\beta_m\not\in T_1(D(m)')^{\perp}$. The sequence
of ideals of the center under consideration takes the form
$$\soc(D(m)') \underbrace{\mbox{~~$=$~~}}_{\mbox{{\small 0}}}
T_1(D(m)')^{\perp}
\underbrace{\mbox{~~$\subset$~~}}_{\mbox{{\small 2}}} Z(D(m)').$$

Secondly, we study the analogous sequence of ideals for the
nonstandard algebra $D(m)$.
\smallskip

{\em Claim.} $\alpha^2\not\in K(D(m))$.
\smallskip

{\em Proof of the Claim.} This can be seen by working out an
explicit basis for the commutator space $K(D(m))$. Let ${\mathcal
B}$ be the monomial basis of $D(m)$ consisting of all pairwise
distinct nonzero paths in the quiver with relations. Set
$\bar{{\mathcal B}}:={\mathcal B}\setminus \{e_1,\ldots,e_m,
\alpha,\alpha^2,s_1,\ldots,s_m\}$, where $e_j$ is the trivial path
at the vertex $j$. Then a $K$-basis of $K(D(m))$ is given as
follows
$$K(D(m))=\langle \bar{{\mathcal B}},\alpha^2-s_2,\ldots,
\alpha^2-s_m, s_m-\alpha^3\rangle_K.$$
Now it is readily checked that $\alpha^2$ can indeed not be written
as a linear combination of these basis elements. $\Diamond$

\smallskip

From the claim we can conclude that $\alpha\not\in T_1(D(m))$.
But the element $\beta_1\ldots\beta_m$ is orthogonal to all
basis elements in the ideal generated by the arrows of the quiver,
except to $\alpha$. In fact,
the product of $\beta_1\ldots\beta_m$ with any path
of length $\ge 1$ becomes 0, except for $\alpha$.
Since $\alpha\not\in T_1(D(m))$, we therefore get that
$$\beta_1\ldots\beta_m\in T_1(D(m))^{\perp}.$$
Hence the sequence of ideal under consideration takes the form
$$\soc(D(m)) \underbrace{\mbox{~~$\subset$~~}}_{\mbox{{\small 1}}}
T_1(D(m))^{\perp} = \langle \beta_1\ldots
\beta_m, \soc(D(m))\rangle_K
\underbrace{\mbox{~~$\subset$~~}}_{\mbox{{\small 1}}} Z(D(m)).$$

Comparing the codimensions of the ideals in these sequences for
$D(m)'$ and $D(m)$, we can finally conclude, by invoking
Proposition \ref{prop:zimmermann}, that the algebras $D(m)'$ and
$D(m)$ are not derived equivalent.
\end{proof}

Summing up, we obtain the following derived equivalence
classification of the representation-finite symmetric algebras.

\begin{theorem}
The algebras $K$, $N_n^{m n}$, $m,n\geq 1$, $D(m)$,
$D(m)^{\prime}$, $m \geq 2$, $T(2,2,r)$, $r \geq 2$, $T(2,3,3)$,
$T(2,3,4)$, $T(2,3,5)$ form a complete set of representatives of
pairwise different derived equivalence classes of connected
representation-finite symmetric algebras.
\end{theorem}


{\sc Ackowledgement.} The first author is grateful to
A.\,Zimmermann for many interesting discussions on the new derived
invariant from \cite{Z}. This started at the Oberwolfach workshop
'Cohomology of finite groups' in September 2005. Many thanks also
to the organizers for the invitation to this event.

\vspace{.2cm}

\begin{footnotesize}
\begin{tabular}{ll}
\normalsize Thorsten \textsc{Holm} \vspace{.1cm}  &
\normalsize Andrzej \textsc{Skowro\'nski} \vspace{.1cm}\\
  Department of Pure Mathematics & Faculty of Mathematics and Computer Science \\
 University of Leeds\hspace{1.2in} &
Nicolaus Copernicus University\\
Leeds LS2 9JT &
Chopina 12/18\\
United Kingdom &
87-100 Toru\'{n}, Poland\\
E-mail: {\tt tholm@maths.leeds.ac.uk} &
E-mail: {\tt skowron@mat.uni.torun.pl} \\
\end{tabular}
\end{footnotesize}

\end{document}